\documentclass[12 pt]{article}  
\newcommand{\hh}{\hat h}
\newcommand{\hcR}{{\cal Q}}

\usepackage{graphicx,amsthm}
\usepackage[usenames]{color}
\usepackage{epsfig}

\newcommand{\good}{effective}
\newcommand{\effective}{fit for control}
\newcommand{\effectiveness}{fitness for control}

\newcommand{\UU}{\mathbf{U}}

\newcommand{\calV}{{\cal V}}
\newcommand{\calW}{{\cal W}}
\newcommand{\calZ}{{\cal Z}}

\newcommand{\fhi}{\varphi}
\newcommand{\Om}{\Omega}
\newcommand{\om}{\omega}

\newcommand{\mt}{\mapsto}
\newcommand{\al}{\alpha}

\newcommand{\Sch}{Schr\"odinger}

\newcommand{\beq}{\begin{equation}}
\newcommand{\eeq}{\end{equation}}
\newcommand{\be}{\begin{equation}}
\newcommand{\ee}{\end{equation}}
\newcommand{\bea}{\begin{eqnarray}}
\newcommand{\eea}{\end{eqnarray}}
\newcommand{\brs}{\begin{eqnarray*}}
\newcommand{\ers}{\end{eqnarray*}}
\newcommand{\ba}{\begin{array}}
\newcommand{\ea}{\end{array}}
\newcommand{\br}{\begin{eqnarray}}
\newcommand{\er}{\end{eqnarray}}

\newcommand{\lp}{\left(}
\newcommand{\rp}{\right)}
\newcommand{\la}{\left\langle}
\newcommand{\ra}{\right\rangle}

\def\EOP{\ \hfill $\Box$}

\renewcommand{\r}[1]{(\ref{#1})}
\def\eps{\varepsilon}
\def\lb{\lambda}
\newcommand{\s}{{\cal S}}

\newtheorem{theorem}{Theorem}[section]

\newtheorem{corol}[theorem]{Corollary}
\newtheorem{lem}[theorem]{Lemma}

\newtheorem{prop}[theorem]{Proposition}
\newtheorem{defn}[theorem]{Definition}

\newtheorem{rem}[theorem]{Remark}
\newcommand{\OO}{{\cal O}}

\newcommand{\N}{{\mathbf N}}
\newcommand{\C}{{\mathbf{C}}}
\newcommand{\R}{{\mathbf{R}}}
\newcommand{\Q}{{\mathbf{Q}}}

\usepackage{graphics,color} 
\usepackage{amsmath} 
\usepackage{amssymb}  

\DeclareMathOperator*{\essinf}{ess\,inf} 
\DeclareMathOperator*{\esssup}{ess\,sup}

\title{\LARGE \bf
Generic controllability properties for the bilinear Schr\"odinger equation\thanks{This work was  supported by the
BQR ``Contr\^ole effectif des syst\`emes quantiques", R\'egion Lorraine -- Nancy Universit\'e.}
}

\author{Paolo Mason\thanks{P. Mason is with
Laboratoire des Signaux et Syst\`emes, Sup\'elec,
3, Rue Joliot Curie, 91192 Gif s/Yvette, France,
        {\tt\small paolo.mason@lss.supelec.fr}},\quad Mario Sigalotti\thanks{M. Sigalotti is with INRIA Nancy - Grand Est, \'Equipe-projet CORIDA, and
Institut \'Elie Cartan, UMR CNRS/INRIA/Nancy Universit\'e, BP 239, 54506 Vand\oe uvre-l\`es-Nancy, France,
        {\tt\small Mario.Sigalotti@inria.fr}}%
}

\begin{document}

\maketitle

\begin{abstract}
In \cite{noi} we proposed a set of sufficient conditions for the
approximate controllability of a discrete-spectrum bilinear
Schr\"odinger equation. These conditions are expressed
in terms of the controlled potential and of the eigenpairs of the
uncontrolled Schr\"odinger operator. The aim of this paper is to
show that these conditions are generic with respect to the
uncontrolled and the controlled potential, denoted respectively by  $V$ and  $W$.
More precisely, we prove that the Schr\"odinger equation is approximately controllable
generically with respect to $W$ when $V$ is fixed and 
also
generically with respect to $V$ when $W$ is fixed and non-constant.
The results are obtained by analytic perturbation arguments and
through the study of asymptotic properties of eigenfunctions.
\end{abstract}

\section{Introduction}
In this paper we consider controlled \Sch\ equations of the type
\be\label{main_EQ} i\frac{\partial \psi}{\partial t}(t,x) = (-\Delta
+V(x)+u(t)W(x))\psi(t,x),\ \ \ \  \ \ \ u(t)\in U,
\ee
where the {\it wave function}
$\psi$ is a map from $[0,+\infty)\times \Omega$ to $\C$ for some domain $\Omega$ of $\R^d$, $d\geq 1$, $V,W$ are 
real-valued functions and $U$ is a nonempty subset of $\R$.
We will assume either that $\Om,V,W$ are bounded and 
$\psi|_{[0,+\infty)\times\partial\Omega}=0$
or that $\Om=\R^d$ and $-\Delta+V+u W$ has discrete spectrum for every $u\in U$.

As proved in \cite{turinici}, the control system
\eqref{main_EQ} is never exactly controllable in $L^2(\Om)$. 
Approximate controllability is known not to hold for some specific system of type 
\r{main_EQ} (see \cite{mirr_rouch}). 
Nevertheless,
several positive controllability results have been proved in recent years.
Among them, let us mention \cite{Beauchard1,beauchard-coron}, where the  exact controllability among regular
enough wave functions is proved for $d=1$, $\Omega$ bounded  and $V=0$, and 
\cite{nersesyan} for the $L^2$-approximate controllability when $\Om$ is bounded under suitable conditions on  $V$ and $W$.
Other controllability results for  related systems have been obtained in \cite{adami_boscain} (when more than one control is available) and in 
\cite{mirrahimi} (when the spectrum has only finitely many discrete eigenvalues). 

In this paper
we focus on the approximate controllability results obtained 
by the authors in collaboration with U.~Boscain and T.~Chambrion 
in \cite{noi}. 
Such results are related to those in \cite{nersesyan} although the sets of sufficient conditions  proposed in the two papers are incomparable and
the techniques used 
are completely different: \cite{nersesyan} applies a control Lyapunov  function approach, whereas \cite{noi} is based on geometric control methods 
for
the Galerkin approximations in the spirit of \cite{Agr-Sar,sergio}. 
  As a consequence, the results in \cite{noi} are valid also in the case in which $\Om=\R^d$ (unlike those in \cite{nersesyan}) 
and when $\Om$ is a   
  manifold and $\Delta$ is the Laplace-Beltrami operator. 
It should also be mentioned that 
the sufficient conditions for approximate controllability proposed in \cite{noi} imply stronger control properties such as 
control of density matrices (see Section~\ref{s-math-f}) or 
tracking of unfeasible trajectories (see \cite{thomas}). 
The aim of this paper is to show that 
such 
sufficient conditions 
are generic.

The genericity issue for the controllability of the \Sch\ equation 
has  already been addressed in the literature. 
In particular, \cite[Lemma 3.12]{nersesyan} proves 
generic $L^2$-approximate controllability with respect to the pair $(V,W)$  
when $d=1$ and $\Om$ is bounded. Newer results can be found in \cite{ner_gen}.
Generic $L^2$-approximate controllability with respect to $(\Om,W)$ in the case $V=0$ is
obtained in \cite{yannick} as a consequence of 
generic properties of the Laplace-Dirichlet operator.  
Other generic controllability results  for a linearized \Sch\ equation can be found in \cite{BCKL} and are further discussed in Section~\ref{concl}.

The difference between our approach and those usually 
adopted to prove generic properties of controlled partial differential equations 
is that, instead of applying local infinitesimal variations, 
we exploit  global, long-range, perturbations. 
The idea is the following: denote by $\Gamma$ the class of systems on which the genericity of a certain property $P$ is studied. 
If we are able to prove the 
existence of at least one  element of $\Gamma$ satisfying $P$, then 
we can {\it propagate} $P$ if some analytic dependence 
properties hold true. In this way we can prove that the property holds in a 
dense subset of $\Gamma$.

The paper is organized as follows: in Section~\ref{s-math-f} we describe the notion of 
solution of \r{main_EQ} (this is a delicate point when $\Om=\R^d$ and $W$ is unbounded) and we recall the approximate controllability results obtained in \cite{noi} (Theorem~\ref{appr-contr}). In particular, we formulate the two conditions ensuring approximate controllability: (i) the 
 spectrum of $-\Delta+V$ is non-resonant and (ii) the potential $W$ couples, directly or indirectly, 
every pair of eigenvectors of $-\Delta+V$. 
We also recall the notion of genericity and we detail the topologies 
with respect to which genericity is considered. 
In Section~\ref{VW} we prove the generic approximate controllability of \r{main_EQ} with respect to the pair $(V,W)$. As an intermediate step, we prove in Proposition~\ref{una_cond} that, generically with respect to $V$, the spectrum of 
$-\Delta+V$ is non-resonant. Section~\ref{un_ingr} refines the results of Section~\ref{VW} by showing that the approximate controllability is generic separately with respect to $V$ or $W$ when $(\Om,W)$ or $(\Om,V)$ is fixed (in the first case, assuming that $W$ is non-constant).  
We conclude the paper with Section~\ref{concl}, where we discuss the 
genericity with respect to $\Om$ of the approximate controllability 
of \r{main_EQ} when $(V,W)$ is fixed. 

\section{Mathematical framework}
\label{s-math-f}

\subsection{Notations and definition of solutions}
Let $\N$ be the set of positive integers. 
For $d\in\N$, denote by $\Xi_d$ the set of
 nonempty, open, bounded and connected subsets of
 $\R^d$ and let  $\Xi_d^\infty=\Xi_d\cup\{\R^d\}$.  
Take $U\subset \R$ and assume that $0$ belongs to $U$.

In the  following we consider the \Sch\ equation \eqref{main_EQ} assuming that
the potentials $V,W$ are taken in $L^\infty(\Omega,\R)$
if $\Om$ belongs to $\Xi_d$ and that 
$V,W\in L^\infty_{\mathrm{loc}}(\R^d,\R)$ and 
$\lim_{\|x\|\to \infty}V(x)+u W(x)=+\infty$ for every $u\in U$ if $\Om=\R^d$. Then, for
every $u\in U$, $-\Delta+V+u W$ (with Dirichlet boundary conditions if $\Om$ is bounded)
 is a skew-adjoint
operator on  $L^2(\Omega,\C)$ with discrete spectrum (see \cite{fried,reed_simon}).
In particular, $-\Delta+V+u W$ generates
a group of unitary transformations
$e^{it(-\Delta+V+u W)}:L^2(\Om,\C)\to L^2(\Om,\C)$.
Therefore,  $e^{it(-\Delta+V+u W)}(\s)=\s$
where  $\s$ denotes the unit sphere of $L^2(\Om,\C)$.

When $\Om$ is bounded, for every  $u\in L^\infty([0,T],U)$
and every $\psi_0\in L^2(\Om,\C)$
there exists a unique weak (and mild) solution $\psi(\cdot;\psi_0,u)\in \mathcal{C}([0,T],L^2(\Om,\C))$.
Moreover, if $\psi_0\in D(A)$ and $u\in \mathcal{C}^1([0,T],U)$ then
$\psi(\cdot;\psi_0,u)$ is differentiable and it is a strong solution of \r{main_EQ}.
(See \cite{BMS} and references therein.)

The situation is more complicated when $\Om=\R^d$ and $W$ is unbounded. However, due to the well-definedness of $e^{it(-\Delta+V+u W)}$ for every $u\in U$ and $t\in \R$, we can always associate a solution 
 \be\label{solu}
\psi(t;\psi_0,u)=e^{(t-\sum_{l=1}^{j-1} t_l)(-\Delta +V+u_j W)}\circ e^{t_{j-1}(-\Delta +V+u_{j-1} W)}\circ \cdots \circ e^{t_1(-\Delta +V+u_1 W)}(\psi_0),
\ee
with every initial condition $\psi_0\in L^2(\Om,\R)$ and every piecewise constant control function $u(\cdot)$. 
Here $\sum_{l=1}^{j-1} t_l\leq t<\sum_{l=1}^{j} t_l$ and 
$$u(\tau)=u_k\quad\mbox{  if ~$\sum_{l=1}^{k-1} t_l\leq \tau<\sum_{l=1}^{k} t_l$}$$
 for $k=1,\dots,j$.

\begin{defn}\label{controllability}
We say that the quadruple $(\Om,V,W,U)$ is approximately controllable
if for every $\psi_0,\psi_1\in \s$ and every $\eps>0$
 there exist $T>0$ and $u:[0,T]\to U$ piecewise constant such that
 $\|\psi_1-\psi(T;\psi_0,u)\|<\eps$.
\end{defn}

It is useful for the applications
to extend the notion of approximate controllability from a single \Sch\ equation to
a (possibly infinite) family of identical systems with different initial conditions, through the study of the evolution of the associated density matrix
(see \cite{albertini,caniggia}).

Let $(\varphi_j)_{j\in \N}$ be an orthonormal basis of $L^2(\Om,\C)$, $(P_j)_{j\in N}$ be a sequence of non-negative numbers such that $\sum_{j=1}^\infty P_j=1$, and denote by $\rho$ the {\it density matrix}
$$
\rho=\sum_{j=1}^{\infty}P_j \varphi_j {\varphi_j}^\ast,
$$
where ${\psi}^\ast(\cdot)=\la\psi,\cdot\ra$, for $\psi\in L^2(\Om,\C)$ and $\la\cdot,\cdot\ra$ denotes the scalar product in $L^2$.
According to the classical definition of density matrix, $\rho$ is a non-negative, self-adjoint operator of trace class (see \cite{reed_simon_1}).
If each $\varphi_j=\varphi_j(t)$ is interpreted as the state of a \Sch\ equation of the form \r{main_EQ}, each equation being characterized by the same potentials $V$ and $W$ and driven by the same piecewise constant control $u=u(t)$, then
the time evolution of the density matrix $\rho=\rho(t)$
is described by
\begin{eqnarray}
\rho(t)=\UU(t) \rho(0) \UU^{\ast}(t)
=\sum_{j=1}^\infty P_j \UU(t) \fhi_j(0) (\UU(t) \fhi_j(0))^\ast
\end{eqnarray}
where the operator $\UU(t)$ is defined by
\be\label{UUU} \UU(t)\psi_0=\psi(t;\psi_0,u)
\ee
and  $\UU^*(t)$ denotes the adjoint of $\UU(t)$.

\begin{defn}
Two density matrices  $\rho_0$ and $\rho_1$ are said to be {\it unitarily equivalent} if there exists a unitary transformation $\UU$ of $L^2(\Om,\C)$ such that  $\rho_1=\UU \rho_0  \UU^\ast$.
\end{defn}
For closed systems  the
problem of connecting two density matrices by a feasible trajectory  makes sense only for pairs  
of  
density matrices
that are  unitarily equivalent.  (The situation is different for open systems,  see for instance \cite{altafini}.)

\begin{defn}\label{densities}
We say that the quadruple $(\Om,V,W,U)$ is approximately controllable
in the sense of density matrices
if for every pair
$\rho_0,\rho_1$ of  unitarily equivalent  density matrices
 and every $\eps>0$
 there exist $T>0$ and  $u:[0,T]\to U$ piecewise constant such that
$\|\rho_1-   \UU(T) \rho_0 \UU(T)^\ast \|<\eps$,
where $\|  \cdot\|$ denotes the operator norm on $L^2(\Om,\C)$ and $\UU$ is defined as in \r{UUU}.
 \end{defn}
It is clear that approximate controllability
in the sense of  density matrices
implies approximate controllability (just take $P_1=1$).

In order to state the approximate controllability result obtained in \cite{noi}, we need to recall
 the following two definitions.
\begin{defn}
The elements of a sequence
$(\mu_{n})_{n\in\N}\subset \R$ are said to be $\Q$-linearly independent (equivalently, the sequence is said to be non-resonant)
if for every $K\in \N$ and $(q_1,\dots,q_K)\in\Q^K\smallsetminus\{0\}$
one has $\sum_{n=1}^K q_n \mu_{n}\ne 0$.
\end{defn}
\begin{defn}
A $n\times n$ matrix $C=(c_{jk})_{1\leq j,k\leq n}$ is said to be connected
 if for every pair of indices $j,k\in\{1,\dots,n\}$ there exists a finite sequence $r_1,\dots,r_l\in\{1,\dots,n\}$ such that $c_{j r_1}c_{r_1 r_2}\cdots c_{r_{l-1}r_l}c_{r_l k}\ne 0$.
\end{defn}
In the following we denote by $\sigma(\Om,V)=(\lb_j(\Om,V))_{j\in\N}$
the non-decreasing sequence of eigenvalues of $-\Delta+V$, counted according to their  multiplicities, and by $(\phi_j(\Om,V))_{j\in\N}$ a corresponding sequence of eigenfunctions. Without loss of generality we can assume that $\phi_j(\Om,V)$ is real-valued for every $j\in\N$. 
Recall moreover that
$(\phi_j(\Om,V))_{j\in\N}$  forms an orthonormal basis
of $L^2(\Om,\C)$. If $j\in\N$ is such that $\lb_j(\Om,V)$ is simple, then 
$\phi_j(\Om,V)$ is uniquely defined up to sign. 

\subsection{Basic facts}

The theorem below recalls the controllability results
obtained 
in \cite[Theorems 3.4, 5.2]{noi}. Here and in the following a map $h:\N\to\N$ is called 
a {\it reordering of $\N$} if it is a bijection.

\begin{theorem}\label{appr-contr}
Let  either (i) $\Omega\in \Xi_d$, $V,W\in L^\infty(\Omega,\R)$ or (ii) $\Om=\R^d$, $V,W\in L^\infty_{\mathrm{loc}}(\R^d,\R)$, $\lim_{|x|\to\infty}V(x)+u W(x)=+\infty$ for every $u\in U$, and 
$|W|$ have at most exponential growth at infinity. 
Assume that 
$U$ contains the interval $[0,\delta)$ for some $\delta>0$,
that the elements of $\big(\lambda_{k+1}(\Om,V)-\lb_k(\Om,V)\big)_{k\in\N}$ are $\Q$-linearly independent and that there 
exists a reordering $h:\N\to \N$ such that
 for infinitely many $n\in \N$ the matrix
$$B^h_{n}(\Om,V,W):=\lp \int_\Omega W(x)\phi_{h(j)}(\Om,V)\phi_{h(k)}(\Om,V)\,dx\rp_{j,k=1}^n$$
 is connected (i.e., $B^h_{n}(\Om,V,W)$ is frequently connected). Then
$(\Omega,V,W,U)$
is approximately controllable
in the sense of  density matrices.
\end{theorem}
\begin{rem}
Notice that, even in the unbounded case, 
 each integral $\int_\Omega W(x)\phi_{j}(\Om,V)\phi_{k}(\Om,V)\,dx$ is well defined. Indeed, when $\Om=\R^d$,  
  the growth of $|W|$ is at most exponential and $e^{a |x|}\phi_j(\R^d,V)\in L^2(\R^d,\R)$ for every $a>0$ and $j\in\N$ (see \cite{agmon}).  
\end{rem}

Let $\calV(\Om)$ be equal to $L^\infty(\Om,\R)$ if $\Om\in\Xi_d$ or
 to 
$\{V\in L^\infty_{\mathrm{loc}}(\R^d,\R)\mid \lim_{\|x\|\to\infty}V(x)=+\infty\}$ if $\Om=\R^d$, and endow $\calV(\Om)$ with the $L^\infty$ topology.

Let us recall
some useful perturbation results
 describing the dependence on $V$ 
 of the spectrum of
the operator $-\Delta +V$.

\begin{theorem}[Continuity]  \label{the-continuity}
Let $\Om\in\Xi_d^\infty$. 
Assume that $\overline{V}$ belongs to $\calV(\Omega)$ and that $\lb_{k}(\Om,\overline{V})$ is
 simple. 
 Then there exists a neighborhood ${\cal N}$ of $\overline{V}$ in $\calV(\Om)$ such that 
 $\lb_{k}(\Om,V)$ is simple for every $V\in {\cal N}$ and 
 $V\mapsto \lb_{k}(\Om,V)$ depends continuously on $V$ on ${\cal N}$. 
 Moreover, 
 the map $V\to \phi_k(\Om,V)$ 
 (defined up to the sign) is continuous  from ${\cal N}$ to $L^2(\Om,\R)$.
 \end{theorem}

 The theorem follows form the remark that, if $V$ tends to $\overline{V}$ in $\calV(\Om)$, then 
 the difference between the two operators $-\Delta+V$ and $-\Delta+\overline{V}$ 
 tends to zero in norm. 
 Therefore, the corresponding resolvents converge in norm, leading to the convergence of eigenvalues and spectral projections (see \cite{katino}). 
 
 We will need in the following a stronger continuity result. 

\begin{prop}\label{the-supercontinuity}
Let $\Om=\R^d$. 
Assume that $\overline{V}$ belongs to $\calV(\R^d)$, $\lb_{k}(\R^d,\overline{V})$ is
 simple, and $W\in L^\infty_{\mathrm{loc}}(\R^d,\R)$ be such that $|W|$ has at most exponential growth. 
 Then there exists a neighborhood ${\cal N}$ of $\overline{V}$ in $\calV(\R^d)$ such that 
 $\lb_{k}(\R^d,V)$ is simple for every $V\in {\cal N}$ and 
 $V\mapsto \sqrt{W}\phi_{k}(\R^d,V)$ (defined up to sign) is a continuous function from ${\cal N}$ to $L^2(\R^d,\C)$. 
 \end{prop}
\proof
Let ${\cal N}$ be a neighborhood of $\overline{V}$
 such that $\lb_{k}(\R^d,V)$ is simple for every $V\in {\cal N}$ (Theorem~\ref{the-continuity}). 
Fix $C,\al>0$ such that $|W(x)|<Ce^{\alpha |x|}$ almost everywhere on $\R^d$. Let, moreover, $\al'$ be a constant larger than $\al$. 

The estimates obtained in \cite{agmon} 
(Theorems 4.1, 4.3 and 4.4)
imply that, up to taking a smaller  ${\cal N}$, there exists $K>0$ such that
$$\int_{\R^d} e^{\al' |x|} \phi_k(\R^d,V)^2 dx< K$$
for every $V\in {\cal N}$.

Let $(V_n)_{n\in\N}$ be a sequence converging to $\overline{V}$ in $\calV(\R^d)$. Since $\al'>\al$, 
given $\eps>0$, 
there exists $R>0$ such that
for every 
$n$ large enough 
\brs
\lefteqn{\int_{\{x\in\R^d\mid \|x\|>R\}} |W(x)| (\phi_k(\R^d,V_n)^2 +\phi_k(\R^d,V)^2)dx <}\\
&& \int_{\{x\in\R^d\mid \|x\|>R\}} C e^{\al |x|}(\phi_k(\R^d,V_n)^2 +\phi_k(\R^d,V)^2)dx
< \eps.
\ers
Moreover, by continuity of $V\mapsto \phi_k(\R^d,V)$ from ${\cal N}$ to $L^2(\R^d,\R)$ and since $W\in L^\infty_{\mathrm{loc}}(\R^d,\R)$, 
$$\int_{\{x\in\R^d\mid \|x\|\leq R\}} |W(x)| (\phi_k(\R^d,V_n)-\phi_k(\R^d,V))^2 dx< \eps$$
for every $n$ large enough.
Therefore, $\sqrt{W} \phi_k(\R^d,V_n)$ converges to $\sqrt{W} \phi_k(\R^d,V)$ 
in $L^2(\R^d,\C)$
as $n$ tends to infinity. 
\EOP

Another crucial  result 
for our needs
concerns analytic perturbation properties.   
\begin{theorem}[{\cite[Chapter VII]{katino}, \cite[Chapter II]{Rellich}}]\label{the-analyticity}
Let $I$ be an interval of $\R$ and $\Om$ belong to $\Xi_d^\infty$.
Assume that $V$ belongs to $\calV(\Om)$ and that 
$\mu\mapsto W_\mu$ is an analytic function from $I$ into $L^\infty(\Omega,\R)$.
Then,
there exist two families of analytic functions
$(\Lambda_k:I \to \mathbf{R})_{k\in\N}$
and $(\Phi_k:I \to L^2(\Omega,\R))_{k\in\N}$    such that
for any $\mu$ in $I$ the sequence  $(\Lambda_k(\mu))_{k\in\N}$ is the 
family of eigenvalues of $-\Delta +V 
+ W_\mu$ counted according to their multiplicities and
$(\Phi_k(\mu))_{k\in\N}$ is an orthonormal basis of
corresponding eigenfunctions.
\end{theorem}
In the following sections we will also make use of the stronger analytic dependence result stated below. 
\begin{prop}
Let $\Om$ belong to $\Xi_d^\infty$ and $\{V_\mu\mid\mu\in [0,1]\}$ be a family of functions in $\calV(\Om)$  such that $V_\mu-V_0$ is analytic in $L^{\infty}(\Om)$ with respect to $\mu$. Let $W\in L^\infty_{\mathrm{loc}}(\Om,\R)$ be such that $|W(x)|\leq C(|V_0(x)|+1)$ for almost every $x\in\Om$, for some positive constant $C$. Then,  if the eigenvalues $\lambda_{j}(\Om,V_\mu)$ and $\lambda_{k}(\Om,V_\mu)$ are simple for $\mu\in (0,1)$, the map
$$
\mu\mapsto \int_\Om W \phi_{j}(\Om,V_\mu)\phi_{k}(\Om,V_\mu)
$$
is analytic from $(0,1)$ to $\R$.
\label{sofisticazzi}
\end{prop}

\proof
First of all notice that, when $\Om$ is bounded, the proposition follows directly from Theorem~\ref{the-analyticity}. Let then $\Om=\R^d$.  
Since the scalar product in $L^2(\R^d,\C)$ is analytic, it is enough to prove that the map $\mu\mapsto\sqrt{W}\phi_k(\R^d,V_\mu)$ is analytic in $L^2(\R^d,\C)$ if $\lambda_k(\R^d,V_\mu)$ is simple for $\mu\in(0,1)$.

Let us first show that, setting $T_\mu=-\Delta+V_\mu$ and 
endowing $D(T_0)$ with the graph norm $\|\phi\|_{T_0}=\|\phi\|_{L^2(\R^d,\C)}+\|T_0\,\phi\|_{L^2(\R^d,\C)}$, the eigenfunctions $\phi_k(\R^d,V_\mu)$ are analytic from $(0,1)$ to $D(T_0)$. 
(This is essentially done in \cite[Theorem~5.6]{crucchi}. We adapt the argument to our case for the sake of completeness.)
Take $\lambda_0$ in the resolvent set of the operator $T_{\mu_0}$, for a fixed $\mu_0\in(0,1)$. For $\mu$ in a neighborhood of $\mu_0$ we have
\begin{equation}
(T_\mu-\lambda_0)^{-1}=(T_{\mu_0}-\lambda_0)^{-1} (\mathrm{Id}+(V_\mu-V_{\mu_0})\,(T_{\mu_0}-\lambda_0)^{-1} )^{-1}\,,\label{compositio} 
\end{equation}
where  $\mathrm{Id}$ denotes the identity. 

Note that $\mu\mapsto (\mathrm{Id}+(V_\mu-V_{\mu_0})\,(T_{\mu_0}-\lambda_0)^{-1} )^{-1}$ is analytic in $\mathcal{L}(L^2(\R^d,\C))$, the space of linear and continuous operators of $L^2(\R^d,\C)$, for $\mu$ in a neighborhood of $\mu_0$.

Notice also that 
$(T_{\mu_0}-\lambda_0)^{-1}$ belongs to $\mathcal L(L^2(\R^d,\C),D(T_0))$, the space of linear and continuous maps from $L^2(\R^d,\C)$ to $D(T_0)$ (endowed with the graph norm), 
as it follows from 
the following series of inequalities:
\begin{eqnarray}
\|(T_{\mu_0}-\lambda_0)^{-1}\phi\|_{T_{0}}
&=&\|(T_{\mu_0}-\lambda_0)^{-1}\phi\|_{L^2(\R^d,\C)}+\|T_0(T_{\mu_0}-\lambda_0)^{-1}\phi\|_{L^2(\R^d,\C)}\nonumber\\
&\leq& \|(T_{\mu_0}-\lambda_0)^{-1}\phi\|_{L^2(\R^d,\C)}+\|\phi+(\lambda_0+V_0-V_{\mu_0})(T_{\mu_0}-\lambda_0)^{-1}\phi\|_{L^2(\R^d,\C)}\nonumber\\
&\leq& \big(\|(T_{\mu_0}-\lambda_0)^{-1}\|+1+\|\lambda_0+V_0-V_{\mu_0}\|_{L^{\infty}(\R^d,\C)}\|(T_{\mu_0}-\lambda_0)^{-1}\|\big)\|\phi\|_{L^2(\R^d,\C)}\,.\nonumber
\end{eqnarray}

Hence $F\mapsto (T_{\mu_0}-\lambda_0)^{-1}F$ is a linear and continuous operator from $\mathcal{L}(L^2(\R^d,\C))$ to $\mathcal L(L^2(\R^d,\C),D(T_0))$. 
It follows from  \eqref{compositio} that $\mu\mapsto(T_\mu-\lambda_0)^{-1}$ is analytic from a neighborhood of $\mu_0$ to  $\mathcal L(L^2(\R^d,\C),D(T_0))$. 

Then 
the eigenfunction $\phi_k(\R^d,V_{\mu})$ is analytic with respect to $\mu$ in $D(T_0)$ 
since
the spectral projection 
$$(2\pi i)^{-1} \oint_{|\lambda-\lambda_k(\R^d,V_{\mu_0})|=\epsilon} (T_{\mu}-\lambda)^{-1}d\lambda\,,$$ 
where 
$\eps$ is small enough, is 
analytic as a function of $\mu$ taking values in $\mathcal L(L^2(\R^d,\C),D(T_0))$. 
(See \cite[Theorem~XII.8]{reed_simon} for details.)

To conclude the proof of the proposition it is enough to check that  the linear map from $D(T_0)$ to $L^2(\R^d,\C)$ sending $\phi$ to $\sqrt{W}\phi$ is continuous, i.e., that there exists $\hat C>0$ such that $\|\sqrt{W}\phi\|_{L^2(\R^d,\C)}\leq \hat C \|\phi\|_{T_0}$ for every $\phi\in D(T_0)$.   
We have 
\begin{eqnarray*}
\|\sqrt{W}\phi\|_{L^2(\R^d,\C)}^2
&\leq& C\left(\|\sqrt{V_0} \,\phi\|_{L^2(\R^d,\C)}^2+\|\phi\|_{L^2(\R^d,\C)}^2\right)\nonumber\\
&\leq& C\left(\big\| \nabla\phi \big\|_{L^2(\R^d,\C)}^2+\big\| \sqrt{V_0}\, \phi \big\|_{L^2(\R^d,\C)}^2+\|\phi\|_{L^2(\R^d,\C)}^2\right)\nonumber\\
&=& C\left(\la (-\Delta+|V_0|) \phi,\phi \ra_{L^2(\R^d,\C)}+\|\phi\|_{L^2(\R^d,\C)}^2\right)\nonumber\\
&\leq&C\left( \|\phi\|_{T_0}\|\phi\|_{L^2(\R^d,\C)}+2\max\{0,-\essinf V_0\}\|\phi\|_{L^2(\R^d,\C)}^2\right)\\
&\leq& \hat C \,\|\phi\|_{T_0}^2,\nonumber
\end{eqnarray*}
where we can take $\hat C=C(1+2\max\{0,-\essinf V_0\})$.  
\EOP

The following proposition 
states the existence of analytic paths of potentials such that the spectrum is simple along them. 

\begin{prop}
Let $\Om$ belong to $\Xi_d^\infty$ and 
$V,Z\in \calV(\Om)$  be such that $Z-V\in L^\infty(\Om,\R)$. 
Then there exists an analytic function $\mu\mapsto W_\mu$  from $[0,1]$ into $L^\infty(\Omega,\R)$ such that $W_0=0$, $W_1=Z-V$ and the spectrum of
$-\Delta+V+W_\mu$ is simple for every $\mu\in(0,1)$.  
\label{teytel-docet}
\end{prop}
\proof
Denote by $\mathcal{C}_0(\Om)$ the subspace of $L^\infty(\Om)$ of continuous real-valued functions vanishing at infinity. Note that $\mathcal{C}_0(\Om)$ is a separable Banach space.
The proof of the proposition is based on \cite[Theorem B]{Teytel}, which guarantees that
the thesis holds true with $W_\mu \in \mathcal{C}_0(\Om)+\R(Z-V)$
provided that, for every $W\in \mathcal{C}_0(\Om)+\R(Z-V)$ and every 
multiple eigenvalue $\lambda$ of $-\Delta+V+W$, there exist two orthonormal eigenfunctions $\phi_1$ and $\phi_2$  pertaining to $\lambda$ such that
the linear functionals 
$$ p\mt \int_\Om p\; (\phi_1^2- \phi_2^2)dx$$
and
$$ p\mt \int_\Om p\, \phi_1\phi_2 dx$$
are linearly independent on $\mathcal{C}_0(\Om)+\R(Z-V)$.

The linear independence of the two functionals can be proved by
taking any pair of orthonormal eigenfunctions  $\phi_1$ and $\phi_2$  pertaining to $\lambda$ and
assuming by contradiction that there exists $(a_1,a_2)\in\R^2\setminus\{(0,0)\}$ such that 
$$ p\mt \int_\Om p (a_1\phi_1^2- a_1\phi_2^2+a_2 \phi_1\phi_2)dx$$
is identically equal to zero on $\mathcal{C}_0(\Om)+\R(Z-V)$.
Hence, $a_1\phi_1^2- a_1\phi_2^2+a_2 \phi_1\phi_2$ must be identically 
equal to zero on $\Om$.

By diagonalizing the quadratic form $a_1 v_1^2- a_1v_2^2+a_2 v_1v_2$ we end up with $c_1,c_2\in \{-1,0,1\}$ and two linearly independent eigenfunctions $\psi_1$ and $\psi_2$ pertaining to $\lambda$ such that $c_1^2+c_2^2>0$ and
$$c_1\psi_1^2+c_2 \psi_2^2\equiv 0.$$

Then for every $x\in \Om$ either $\psi_1(x)=\psi_2(x)$ or $\psi_1(x)=-\psi_2(x)$.
Thanks to the unique continuation property (see \cite{simon_semigroups}) we have $\psi_1=\pm \psi_2$, contradicting the linear independence of $\psi_1$ and $\psi_2$. 
\EOP

\subsection{Genericity: topologies and definitions}

From now on we will write simply $L^p(\Om)$ or $L^p_{\mathrm{loc}}(\Om)$
to denote $L^p(\Om,\R)$ or $L^p_{\mathrm{loc}}(\Om,\R)$ respectively.
Similarly, $H^p(\Om)$ and $H^p_0(\Om)$ will denote $H^p(\Om,\R)$ and $H^p_0(\Om,\R)$.

For every $\Om\in\Xi_d^\infty$ 
let 
$\calW(\Om)$ 
 be equal to $L^\infty(\Om)$ if $\Om\in\Xi_d$ or to 
$$\{W\in L^\infty_{\mathrm{loc}}(\R^d)\mid \esssup_{x\in\R^d}\frac{\log(|W(x)|+1)}{\|x\|+1}<\infty\}$$ 
if $\Om=\R^d$.
In both cases endow $\calW(\Om)$  with the $L^\infty$ topology. 
Denote 
$$\calZ(\Om,U)=\{(V,W)\mid V\in\calV(\Om),\ W\in\calW(\Om), V+uW\in{\cal V}(\Om)\mbox{ for every $u\in U$}\}$$
and endow $\calZ(\Om,U)$ with the product $L^\infty$ topology. 
We also introduce, for every $V\in\calV(\Om)$ and every $W\in\calW(\Om)$,  the topological subspaces of $\calV(\Om)$ and $\calW(\Om)$ defined,
with a slight abuse of notation, by
\brs 
\calV(\Om,W,U)&=&\{\tilde V\in\calV(\Om)\mid (\tilde V,W)\in\calZ(\Om,U)\},\\
\calW(\Om,V,U)&=&\{\tilde W\in\calW(\Om)\mid (V,\tilde W)\in\calZ(\Om,U)\}.
\ers
Notice that neither $\calV(\Om,W,U)$ nor $\calW(\Om,V,U)$ is empty.
Moreover, both $\calV(\Om,W,U)$ and $\calW(\Om,V,U)$ are invariant by the set addition with $L^\infty(\Om)$. 
In particular, they are open in $\calV(\Om)$ and $\calW(\Om)$ respectively and 
they coincide with $L^\infty(\Om)$ when $\Om\in\Xi_d$. 

Theorem \ref{appr-contr} motivates the following definition.
\begin{defn}
Let $V\in\calV(\Om)$ and $W\in \calW(\Om)$. We say that $(\Om,V,W)$ is
{\it \effective} if $(\lb_{k+1}(\Om,V)-\lb_k(\Om,V))_{k\in\N}$ 
 is non-resonant and
$B^h_{n}(\Om,V,W)$ is frequently connected for some reordering $h$.
Let $(V,W)$ be an element of $\calZ(\Om,U)$. We say that the quadruple $(\Om,V,W,U)$ is {\it \good} if
$(\Om,V+uW,W)$ is \effective\ for some $u$ such that $[u,u+\delta)\subset U$ for some $\delta>0$. 
\label{sake}
\end{defn}
Theorem \ref{appr-contr} can then
be rephrased by 
saying that
being \good\ is a sufficient condition for
controllability in the sense of the density matrices. The rest of the paper deals with the genericity of the notions introduced in Definition~\ref{sake}.

Let us recall that a topological space $X$ is called a {\it Baire space} 
if the intersection of countably many open and dense subsets of $X$  is dense in $X$.
Every complete metric space is a Baire space. 
The intersection of countably many open and dense subsets of a Baire space is called a
{\it residual subset} of $X$. 
Given a Baire space $X$, 
a boolean function $P:X\to \{0,1\}$ is said to be a {\it generic property} if
there exists a residual subset $Y$ of $X$ such that every $x$ in $Y$ satisfies property $P$, that is, $P(x)=1$.

\section{
The triple $(\Om,V,W)$ is generically \effective\
with respect to the pair $(V,W)$}\label{VW}

Let us start by recalling a known result on the generic simplicity of eigenvalues for \Sch\ operators on bounded domains (see \cite{Aalbert}).
\begin{prop}[Albert]\label{albert}
Let $\Omega$ belong to $\Xi_d$.
For every $k\in \N$  the set
$$
\{V\in  L^\infty(\Omega)\mid \lb_1(\Om,V),\dots,\lb_{k}(\Om,V)\mbox{ are simple}\}
$$
is open and dense in $L^\infty(\Omega)$.
Hence, the spectrum $\sigma(\Om,V)$ is simple generically with respect
to $V$.
\end{prop}

For every $\Om\in\Xi_d^\infty$ and every $k\in\N$, let
\be\label{Rk}
{\cal R}_k(\Om)=\{V\in  \calV(\Omega)\mid \lb_1(\Om,V),\dots,\lb_{k}(\Om,V)\mbox{ are simple}\}.
\ee
We generalize Proposition~\ref{albert}
as follows.
\begin{prop}\label{una_cond}
Let $\Omega$ belong to $\Xi_d^\infty$. 
For every $K\in \N$ and
$q=(q_1,\dots,q_K)\in\Q^K\setminus\{0\}$, the set 
\be\label{Oq}
\OO_q(\Om)=\left\{V\in  {\cal R}_K(\Om)
\mid \sum_{j=1}^{K}q_j\lb_j(\Om,V)\ne
0\right\} 
\ee 
is open and dense in $\calV(\Omega)$. Hence, the
spectrum $\sigma(\Om,V)$ forms a
non-resonant family generically with respect to $V$.
\end{prop}

The proof of Proposition~\ref{una_cond} is based on the following lemma.
\begin{lem}\label{buco}
Let $\Omega$ belong to $\Xi_d^\infty$
and $\omega$ be a nonempty,  open
subset compactly contained in $\Omega$ and whose
boundary is Lipschitz continuous.
Let 
$v\in L^\infty(\om)$ and $(V_k)_{k\in\N}$ be a sequence in 
$\calV(\Om)$
such that $V_k|_\om\to v$ in $L^\infty(\omega)$ as $k\to\infty$  and
$\lim_{k\to\infty}\inf_{\Omega\setminus\omega}V_k=+\infty$.
Then, for every $j\in\N$,  $\lim_{k\to\infty}\lb_j(\Om,V_k)=\lb_j(\om,v)$.
Moreover, if $\lb_j(\om,v)$ is simple then (up to the sign) $\phi_j(\Om,V_k)$ and $\sqrt{V_k}\phi_j(\Om,V_k)$ converge respectively to $\phi_j(\om,v)$ and $\sqrt{v}\phi_j(\om,v)$ in $L^2(\Om,\C)$ as $k$ goes to infinity, 
where $\phi_j(\om,v)$ is identified with its extension by zero outside $\omega$. 
\end{lem}
\noindent {\it Proof of Lemma~\ref{buco}.}
Without loss of generality
we can assume that $v\geq 1$ on $\om$
and $V_k\geq 1$ on $\Om$.
Indeed, for $k$ large enough  $V_k\geq 1$ on $\Om\setminus\om$ 
and, for what concerns the values of $V_k$ and $v$ on $\om$, it suffices to notice that 
if we replace
$v$ by
$v+c$   and each $V_k$
by
$V_k+c$ 
with 
$$c=\max\{\|v\|_{L^\infty(\om)},\|V_k\|_{L^\infty(\om)}\mid k\in \N\}+1,$$
then we simply operate a shift of the spectra $\sigma(v,\om)$ and $\sigma(V_k,\Om)$ by the common 
constant $c$.  
Therefore, 
the operators
$-\Delta+v$ 
and 
$-\Delta+V_k$ 
are invertible and their inverses, denoted respectively by $a:L^2(\om)\to L^2(\om)$ and $A_k:L^2(\Om)\to L^2(\Om)$, are compact and have uniformly bounded norm.

Let $A:L^2(\Om)\to L^2(\Om)$ be the operator
associating with
$f\in L^2(\Om)$ the extension by zero
on $\Om\setminus  \om$ of $a(f|_\om)$.
Let us prove that $A_k$ converges pointwise
to $A$.

Fix $f\in L^2(\Omega)$
and denote, for every $k\in\N$,
$w_k=A_k f$ and $w=Af$.
Then $w_k\in H^1_0(\Omega)$ 
is the weak solution of
$$
(-\Delta+V_k) w_k =f\ \ \ \ \ \mbox{ in }\Om
$$
and $w|_\om\in H^1_0(\om)$ is the  weak solution of
\be\label{uniq}
(-\Delta+v) w =f\ \ \ \ \ \mbox{ in }\om.
\ee
We must prove that
 $\|w_k-w\|_{L^2(\Om)}$
 tends to zero as $k$ goes to infinity.

By definition of weak solution, for every $\fhi\in H^1_0(\Om)$,
\be\label{w1}
\int_\Om \nabla w_k\cdot\nabla\fhi+\int_\Om V_k w_k \fhi=\int_\Omega f \fhi
\ee
and, similarly, for every $\psi\in H^1_0(\om)$,
$$
\int_\om \nabla w\cdot\nabla\psi+\int_\om v w \psi=\int_\omega f \psi.
$$

Taking $w_k$ as $\fhi$ in \r{w1} we easily get that the sequence $w_k$ is uniformly bounded in $H^1_0(\Om)$. Denote by $w^*$
the limit of a subsequence of $w_k$ weakly converging  in $H^1_0(\Om)$ (whose existence is guaranteed by Banach-Alaoglu theorem). 
With a slight abuse of notation, let us identify $w_k$ with its weakly converging subsequence.

The definition of weak convergence and \r{w1}
imply that,
for every $\psi\in H^1_0(\om)$,
$$
\int_\om \nabla w^*\cdot\nabla\psi+\int_\om v w^* \psi=\int_\omega f \psi.
$$
Taking again $\fhi=w_k$ in \r{w1}, we notice that
$\{\sqrt{\inf_{\Omega\setminus\omega} V_k}\|w_k\|_{L^2(\Om\setminus\om)}\}_{k\in\N}$ is a bounded sequence in $\R$. Hence, $w_k\to 0$ in $L^2(\Om\setminus\om)$ 
and thus $w^*=0$ on $\Om\setminus\om$.
Recall
that,
since the boundary of $\om$ is Lipschitz continuous, then any 
$H^1$ function which is defined in a neighborhood of $\om$ and which annihilates outside $\om$ belongs to $H^1_0(\om)$ (see \cite[Lemma 3.2.15]{henrot-pierre}).
Hence $w^*\in H^1_0(\om)$ 
coincides with $w$, the unique weak solution of \r{uniq}.
Since a $H^1$-weakly converging sequence converges $L^2$-strongly on bounded sets 
(see, for instance, \cite[Theorem 8.6]{lieb_loss}), then we have proved the convergence of $w_k$ to $w$ in $L^2(\Om)$, that is, the pointwise convergence of $A_k$ to $A$.

We claim that the family of operators $\{A_k\mid k\in\N\}$ is collectively compact. Recall that  
 $\{A_k \mid k\in\N\}$ is collectively compact if for every $(f_k)_{k\in\N}$ in the unit ball of $L^2(\Om)$, the set $\{A_k f_k\mid k\in\N\}$ is pre-compact (see \cite{anselone}). The proof of this fact is quite classical and the argument proposed here follows a similar one given in \cite[Lemma 5.1]{allaire}.
Let $z_k=A_k f_k$ and notice that
\be\label{liboni}
 \|\nabla z_k\|^2_{L^2(\Om)}+\|\sqrt{V_k}z_k\|^2_{L^2(\Om)}=\int_\Om f_k z_k.
 \ee
Since the $A_k$'s are uniformly bounded, then 
the right-hand side of \r{liboni} is uniformly bounded. 
Thus, up to extracting a subsequence, $z_k$ weakly converges to some $z$ in $H^1_0(\Om)$. 
If $\Om\in\Xi_d$ then $z_k\to z$  
strongly in $L^2(\Om)$ which proves the collective compactness of $A_k$ in this case. 
When $\Om=\R^d$ let $\overline{V}(x)=\inf_{k\in\N}V_k(x)$ for every $x\in\R^d$. 
Then $\overline{V}$ belongs to $L^\infty_{\mathrm{loc}}(\R^d)$, $\overline{V}\geq 1$ almost everywhere and 
$\lim_{\|x\|\to \infty} \overline{V}(x)=+\infty$. 
In order to prove this last property, assume by contradiction that there exists a sequence $x_k$  such that 
$\lim_{k\to\infty}\|x_k\|=\infty$ and $\overline{V}(x_k)$ is uniformly bounded. Then there exists 
a subsequence $x_{k_j}$ such that either 
$v_j=V_m(x_{k_j})$ is uniformly bounded for some  $m\in\N$   or $V_{m_j}(x_{k_j})$ is uniformly bounded for some unbounded sequence $(m_j)_{j\in\N}$ in $\N$.
 The contradiction follows in the first case from the fact that $V_m\in\calV(\R^d)$, 
 while in the second case it is a consequence of  the convergence of $\inf_{\R^d\setminus\om}V_k$ to infinity as $k$ goes to infinity. 
Then, for every $\rho>0$, 
\be\label{guaina}
\|z_k-z\|^2_{L^2(\R^d)}\leq \int_{\{\overline{V}< \rho\}}(z_k-z)^2+\frac{1}{\rho}\int_{\{\overline{V}\geq \rho\}}(z_k-z)^2\,\overline{V}
\ee
 and $\{\overline{V}< \rho\}$ is bounded. 
 It follows from \r{liboni} that the $L^2$-norm of $z_k\sqrt{\overline{V}}$ on $\R^d$ is uniformly bounded with respect to  $k$. Since 
 $z_k$ converges $H^1$-weakly  to $z$ in $\R^d$, and therefore $L^2$-strongly on each compact set, then 
$z \sqrt{\overline{V}}$ belongs to $L^2(\R^d)$.  
We deduce that $\int_{\{\overline{V}\geq \rho\}}(z_k-z)^2\,\overline{V}$ is uniformly bounded with respect to $k$ and thus, for $\rho$ large enough, $(1/\rho)\int_{\{\overline{V}\geq \rho\}}(z_k-z)^2\,\overline{V}$ is arbitrarily small, uniformly with respect to $k$.
Since $\{\overline{V}< \rho\}$ is bounded, then, for any fixed $\rho$, 
$z_k\to z$ in $L^2(\{\overline{V}< \rho\})$. 
It follows from \r{guaina} that $z_k$ converges to $z$ in $L^2(\R^d)$, concluding the proof of the 
collective compactness of $\{A_k\}_{k\in\N}$.

Theorems 4.8 and 4.11 in \cite{anselone} guarantee that 
$\lb_j(\Om,V_k)$ converges to $\lb_j(\om,v)$ as $k$ goes to infinity for every $j\in\N$ and 
that, if $\lb_j(\om,v)$ is simple, then (up to the sign)
$\lim_{k\to\infty}\phi_j(\Om,V_k)=\phi_j(\om,v)$ in $L^2(\Om)$.

To conclude the proof 
we observe that  it is enough to show that the restriction of $\sqrt{V_k}\phi_j(\Om,V_k)$ to $\Om\setminus\om$ converges to zero in $L^2(\Om\setminus \om)$ as $k$ goes to infinity, where $j$ is such that $\lb_j(\om,v)$ is simple.

Since $\phi_j(\Om,V_k)$ satisfies
$$(-\Delta+V_k)\phi_j(\Om,V_k) =\lb_j(\Om,V_k)\phi_j(\Om,V_k)$$
in the weak sense, then taking
$\phi_j(\Om,V_k)$ as test function we have
\be\label{tutto-torna}
\int_\Om \|\nabla \phi_j(\Om,V_k)\|^2
+\int_\Om V_k \phi_j(\Om,V_k)^2 =\lb_j(\Om,V_k).
\ee
In particular the sequence $(\phi_j(\Om,V_k))_{k\in\N}$ is uniformly bounded in $H_0^1(\Om)$ and so 
$\phi_j(\Om,V_k)$ converges to 
$\phi_j(\om,v)$ not only strongly in $L^2(\Om)$ 
but also weakly in  $H_0^1(\Om)$.
Hence,
\brs
\limsup_{k\to\infty} \|\sqrt{V_k} \phi_j(\Om,V_k)\|^2_{L^2(\Om\setminus \om)}&=&\lb_j(\om,v)-\int_\om v \phi_j(\om,v)^2-\liminf_{k\to\infty} \|\nabla \phi_j(\Om,V_k)\|^2_{L^2(\Om)}\\
&\leq&\lb_j(\om,v)-\int_\om v \phi_j(\om,v)^2- \|\nabla \phi_j(\om,v)\|^2_{L^2(\Om)}.
\ers
The last term of the above inequality is equal to zero, as it follows 
from the analogous of \r{tutto-torna} for $\phi_j(\om,v)$. 
\EOP

\noindent {\it Proof of Proposition~\ref{una_cond}}.
The second part of the statement clearly follows from the first one, since
\[\{V\in \calV(\Om)\mid \sigma(\Om,V) \mbox{ non-resonant}\}=
\cap _{q\in \cup_{K\in\N}(\Q^K\setminus\{0\})}\OO_q(\Om).
\]

For each 
$K\in \N$ and $q=(q_1,\dots,q_K)\in\Q^K\setminus\{0\}$,
the openness of $\OO_q(\Om)$ in $\calV(\Om)$
follows directly from the continuity of the eigenvalues
on $V$. (See Theorem~\ref{the-continuity}.)

We prove the  density of $\OO_q(\Om)$ in $\calV(\Om)$ by an analytic perturbation argument.
Fix $V\in \calV(\Om)$.
Let $\om$ be a $d$-orthotope compactly contained in $\Om$ and $v$ a measurable bounded function on $\om$ such that $\sigma(\om,v)$ is  non-resonant. 
(The existence of such $\om$ and $v$ is obtained in \cite[Section 6.3]{noi} for $d=3$ and the proof extends with no extra difficulty to the general case $d\in \N$.)

Let us consider a sequence $(V_k)_{k\in\N}\in\mathcal{V}(\Omega)$
such that $V_k-V\in L^\infty(\Omega),\ \forall k\in\N$ and such that
$V_k|_\om\to v$ in $L^\infty(\omega)$ as $k\to\infty$  and
$\lim_{k\to\infty}\inf_{\Omega\setminus\omega}V_k=+\infty$. By
Lemma~\ref{buco} we have that $\lim_{k\to\infty}\sum_{j=1}^{K}q_j\lb_j(\Om,V_k)=\sum_{j=1}^{K}q_j\lb_j(\omega,v)\neq 0$ so that $V_{\bar k}\in\mathcal{O}_q$ for some $\bar k$  large enough. 

By Proposition~\ref{teytel-docet} we can construct an analytic path  $\mu\mapsto W_\mu$  from $[0,1]$ into $L^\infty(\Omega)$
such that $W_0=0$, $W_1=V_{\bar k}-V$ and the spectrum of $-\Delta+V+W_\mu$
is simple for every $\mu\in(0,1)$. This, together with Theorem~\ref{the-analyticity}, implies that the map $\mu\mapsto \sum_{j=1}^{K}q_j\lb_j(\Om,V+W_\mu)$, which is different from zero at $\mu=1$, is analytic and thus different from zero almost  everywhere. Hence, $\mathcal{O}_q$ is dense in $\mathcal{V}(\Omega)$.
\EOP

The following theorem
extends the analysis from $V$ to the pair $(V,W)$,
combining the
generic non-resonance of the spectrum of $-\Delta+V$ with the genericity of the  connectedness 
of the matrices $B^h_{n}(\Om,V,W)$.
\begin{theorem}\label{wrtVW}
Let $\Om$ belong to $\Xi_d^\infty$.
Then, generically with respect to $(V,W)\in \calZ(\Omega,U)$
the triple $(\Om,V,W)$ is \effective.
\end{theorem}
\proof
We proved in Proposition~\ref{una_cond} that each ${\cal R}_k(\Om)$, defined in \r{Rk},
is open and dense in $\calV(\Om)$.
If $V$ belongs to ${\cal R}_k(\Om)$, then the eigenfunctions $\phi_1(\Om,V),\dots,\phi_k(\Om,V)$ are uniquely defined in $L^2(\Om)$ 
up to sign.
It makes sense, therefore, to define
\brs
{\cal U}_k(\Om,U)=\{ (V,W)\in \calZ(\Om,U)&\mid&
V\in{\cal R}_k(\Om),
\\
&& \int_\Om W \phi_{j_1}(\Om,V)\phi_{j_2}(\Om,V)\ne 0\mbox{ for every }1\leq j_1,j_2\leq k\}.
\ers

Let $1\leq j_1,j_2\leq k$.
As it follows from the unique continuation property (see \cite{simon_semigroups}),  the product $\phi_{j_1}(\Om,V)\phi_{j_2}(\Om,V)$ is a nonzero function on $\Om$.  
The set of potentials $W$ belonging to $\calW(\Om)$ that are not orthogonal to $\phi_{j_1}(\Om,V)\phi_{j_2}(\Om,V)$ is therefore 
open and dense in $\calW(\Om)$. Intersecting all such sets as $j_1$ and $j_2$ vary in $\{1,\dots,k\}$ we obtain
again an open and dense subset of $\calW(\Om)$. 
Hence,  ${\cal U}_k(\Om,U)$
is dense in $\calZ(\Om,U)$.
Its openness, moreover, follows from 
Proposition~\ref{the-supercontinuity}. 

The proof is concluded by noticing that
$(\Om,V,W)$ is \effective\ if $(V,W)$ belongs to
$$\lp \cap_{k\in\N}{\cal U}_k(\Om,U)\rp \cap \lp \cap_{q\in \cup_{K\in \N}\Q^K\setminus \{0\}}\{(V,W)\in \calZ(\Om,U)\mid V\in \OO_q(\Om)\}\rp,$$
which is a countable intersection of open and dense subsets of $\calZ(\Om,U)$.
\EOP

\section{Generic controllability with respect to one single
argument}\label{un_ingr}

The following technical result will 
play a crucial role 
in the later discussion.

\begin{lem}\label{AC}
Let $\Om$ belong to $\Xi_d^\infty$ and $Z$ be a non-constant absolutely continuous
function on $\Om$.
Then there exist
$\om\in \Xi_d$ compactly contained in $\Om$ with Lipschitz continuous boundary 
 and
a reordering $h:\N\to\N$ such that
$\sigma(\om,0)$ is simple and 
\be\label{eq_snake}
\int_\om Z\phi_{h(l)}(\om,0)\phi_{h(l+1)}(\om,0)\ne 0
\ee
for every $l\in \N$.
\end{lem}
\proof
Let $\bar x\in \Om$ be such that $\nabla Z(\bar x)$ exists and is different from zero.
Up to a change of coordinates  $\bar x=0$ and each component of $\nabla Z(0)=(\partial_1 Z(0),\dots,\partial_d Z(0))$ is different from zero.

Take as
$\om$  an orthotope of the type $(0,r_1)\times\cdots\times (0,r_d)$ such that
$\sigma(\om, 0)$ is simple.
This is true, for instance, if 
$$\prod_{{\scriptsize \ba{c}1\leq i\leq d\\
i\ne j\ea}}r_i^2,\quad\quad j=1,\dots,d,$$
are $\Q$-linearly independent.
Let $r=(r_1,\dots,r_d)$. 

The choice of $r$ guarantees that $\sigma(\al \om,0)$ is simple 
for every $\al>0$. Therefore, the eigenfunctions of $-\Delta$ on $\al\om$ are
uniquely defined up to sign by 
$$\psi^\al_{k}(x_1,\dots,x_d)=\frac{2^{\frac d2}}{\al^{\frac d2}\sqrt{\prod_{i=1}^d r_i}}\prod_{i=1}^d \sin\lp\frac{k_i x_i \pi}{\al r_i}\rp$$
where $k=(k_1,\dots,k_d)$ belongs to $\N^d$. 

Denote by $f_i(k)$ the element of $\N^d$ obtained from $k$ by adding $1$ to its $i^\mathrm{th}$ component.
By construction $Z(x)=Z(0)+\sum_{j=1}^d x_j \partial_j Z(0)+z(x)$ with $\lim_{x\to0}z(x)/\|x\|=0$.  
Hence, 
\[
\int_{\al \om} Z(x) \psi^\al_{k}(x)\psi^\al_{f_i(k)}(x)dx= \frac{2\partial_iZ(0)}{{\al r_i}}\int_0^{\al r_i} x_i \sin\lp  \frac{k_i x_i \pi}{\al r_i}\rp\sin\lp  \frac{(k_i+1) x_i \pi}{\al r_i}\rp dx_i+\fhi_{k,i}(\al)
\]
with $|\fhi_{k,i}(\al)|\leq \fhi(\al)$ and $\fhi$, independent of $k$ and $i$, satisfies $\lim_{\al\to0}\fhi(\al)/\al=0$. (One can take, for instance, $\fhi(\al)=\sup_{x \in\al\om}|z(x)|$.)

Notice now that
\[\frac{1}{{\al^2 r_i^2}}\int_0^{\al r_i} x_i \sin\lp  \frac{k_i x_i \pi}{\al r_i}\rp\sin\lp  \frac{(k_i+1) x_i \pi}{\al r_i}\rp dx_i=-\frac{4 k_i(1+k_i)}{(1+2 k_i)^2\pi^2}
\]
and that
\[
\lim_{k_i\to \infty}-\frac{4 k_i(1+k_i)}{(1+2 k_i)^2\pi^2}=-\frac{1}{\pi^2}.
\]
Therefore, for $\al$ small enough, we have
 $$\int_{\al\om} Z(x) \psi^\al_{k}(x)\psi^\al_{f_i(k)}(x)dx\ne 0$$
for every $k\in\N^d$ and every $i\in\{1,\dots,d\}$.

\begin{figure}\begin{center}\input{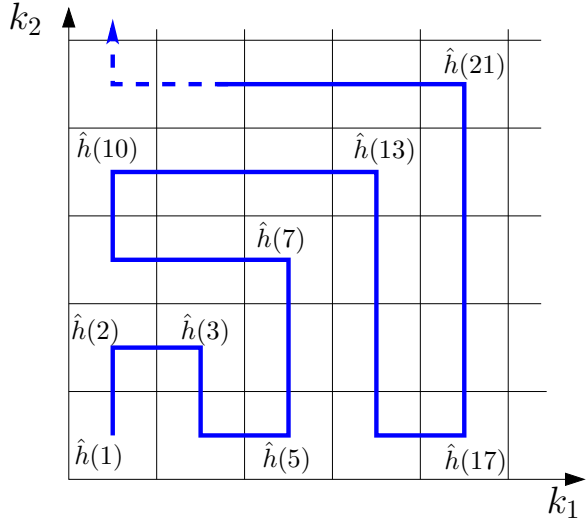}\caption{A possible choice of the function $\hh$ for $d=2$.}\label{fig-snake}\end{center}\end{figure}

We are left to prove that 
there exists a
bijection $\hh:\N\to \N^d$ such that 
$\|\hh(j+1)-\hh(j)\| = 1$ for every $j\in \N$.
This can be interpreted by saying 
that an infinite-length snake as in \cite{snake}  
can fill $\N^d$ (see Figure~\ref{fig-snake}).


 We claim that the following holds:
 \begin{quote}
 For every $m=(m_1,\dots,m_d)\in\N^d$ such that each $m_j$ is odd, 
there exists a bijection 
$$\hh_m:\left\{1,\dots,\prod_{i=1}^d m_i\right\}\longrightarrow \{1,\dots,m_1\}\times \cdots\times\{1,\dots,m_d\}$$
such that  $\|\hh_m(j+1)-\hh_m(j)\|=1$ for $j=1,\dots,\prod_{i=1}^d m_i-1$, $\hh_m(1)=(1,\dots,1)$ and $\hh_m(\prod_{i=1}^d m_i)=m$. 
Moreover,  
if we define $\bar m=(\bar m_1,\dots,\bar m_d)=(m_1,\dots,m_{p-1},m_p+2,m_{p+1},\dots,m_d)$ for some $p\in\{1,\dots,d\}$, then the map $\hh_m$ can be extended to a bijection  
$$\hh_{\bar m}:\left\{1,\dots,\prod_{i=1}^d{\bar m}_i\right\}\longrightarrow \{1,\dots,{\bar m}_1\}\times \cdots\times\{1,\dots,{\bar m}_d\}$$
verifying  $\|\hh_{\bar m}(j+1)-\hh_{\bar m}(j)\|=1$ for $i=1,\dots,\prod_{i=1}^d{\bar m}_i -1
$, $\hh_{\bar m}(1)=(1,\dots,1)$ and $\hh_{\bar m}(\prod_{i=1}^d{\bar m}_i)={\bar m}$. 
\end{quote}

To prove the existence of such $\hh_m$ we proceed by induction. For $d=1$ the claim is trivial. Let now $d=\bar d>1$ and assume that the claim is true for $d=\bar d-1$. Let $m=(m_1,\dots,m_{\bar d})\in\N^d$ with $m_i$ odd for every $i=1,\dots,{\bar d}$. The first part of the claim on $\hh_m$ is obvious when $m_{\bar d}=1$ by the inductive assumption. If $m_{\bar d}>1$ we consider a function $\hh_{\hat m}$ satisfying the first part of the claim with $\hat m=(m_1,\dots, m_{\bar d-1})$. 
For simplicity denote $\mu=\prod_{l=1}^{\bar d-1}m_l$.
Then the map
$$\hh_m(i)=\left\{\ba{ll} \!\!\!
(\hh_{\hat m}(i-j\,\mu),j+1)&\mbox{for } i=j\,
\mu+1,\dots,(j+1)\,\mu,\, j\mbox{ even,}\\
\!\!\!(\hh_{\hat m}((j+1)\,\mu-i+1),j+1)&\mbox{for } i=j\,\mu+1,\dots,(j+1)\,\mu,\, j\mbox{ odd,}
\ea\right.$$
satisfies the required properties. 
As for the second part of the claim on $\hh_m$,
let us assume without loss of generality that $p=\bar d$ (if this is not the case it is enough to reorder the indices of $m_1,\dots, m_{\bar d}$). Then the map
$$\hh_{\bar m}(i)=\left\{\ba{ll} 
\hh_m(i)&\mbox{for }1\leq i\leq \mu m_{\bar d} \\
(\hh_{\hat m}(\mu(m_{\bar d}+1)-i+1), m_{\bar d}+1)&\mbox{for } \mu m_{\bar d}+1\leq i\leq \mu (m_{\bar d}+1)\\
(\hh_{\hat m}(i-\mu(m_{\bar d}+1)),m_{\bar d}+2)&\mbox{for } \mu(m_{\bar d}+1)+1\leq i\leq \mu(m_{\bar d}+2)
\ea\right.$$
satisfies the required properties.

To conclude the proof of the existence of $\hh$ it is enough to consider a sequence of $d$-uples $m(l)=(m_1(l),\dots,m_d(l))$ with positive odd components, such that for every $l$ there exists $p$ with $m_p(l+1)=m_p(l)+2$ and $m_j(l+1)=m_j(l)$ for $j\neq p$, and moreover $m_j(l)$ goes to infinity as $l$ goes to infinity for any fixed $j$. (Take, for instance, $p=l(\mathrm{mod}~d)+1$.) The map $\hh$ is then obtained by extending inductively each map $\hh_{m(l)}$ to a map $\hh_{m(l+1)}$.
\EOP

\subsection{The triple $(\Om,V,W)$ is generically \effective\
with respect to $V$
}

Let $\Om\in\Xi_d^\infty$ and fix $W\in\calW(\Om)$. Let us consider the following subspace of $\calV(\Om)$ 
$$\hat \calV(\Om,W)=\{V\in\calV(\Om)\mid \esssup_{x\in\Om} \frac{|W(x)|}{|V(x)|+1}<+\infty\}\,.$$  

\begin{theorem}\label{V-s}
Let $\Omega$ belong to $\Xi_d^\infty$ and $W\in\calW(\Om)$ be non-constant and absolutely continuous. 
Then, generically with respect to
$V$ in 
$\hat \calV(\Om,W)$,
the triple 
$(\Om,V,W)$
is \effective.
\end{theorem}
\proof
We will denote by  $\hcR_n(\Om,W)$ the set
of potentials $V\in\hat  \calV(\Om,W)$ such that  for  every pair of indices $j,k\in\{1,\dots,n\}$ there exists a finite sequence $r_1,\dots,r_l\in\N$ such that 
$r_1=j$, $r_l=k$, $\lb_{r_i}(\Om,V)$ is simple for every $i=1,\dots,l$, and 
$$\int_\Om W \phi_{r_i}(\Om,V)\phi_{r_{i+1}}(\Om,V)\ne 0$$
for every $i=1,\dots,l-1$.

The openness of $\hcR_n(\Om,W)$ follows from
Proposition~\ref{the-supercontinuity}.

As for its density, apply Lemma~\ref{AC} with $W$ playing the role of  $Z$. Then there exist $\om\in\Xi_d$ with Lipschitz boundary and compactly contained in
$\Om$, and a reordering $h$ of $\N$ such that $\sigma(\om,0)$ is simple and
\begin{equation}
\int_\om W \phi_{h(l)}(\om,0)\phi_{h(l+1)}(\om,0)\ne 0
\label{eq-disug}
\end{equation}
for every $l\in \N$.

Given $\overline{V}\in\hat\calV(\Om,W)$,
let $(V_k)_{k\in\N}$ be the sequence in $\calV(\Om)$ defined by $V_k=0$ in $\omega$ and $V_k=\overline V+k$ in $\Om\setminus\omega$.

Since we know from Lemma~\ref{buco} that $\|\sqrt{V_k}\phi_j(\Om,V_k)\|_{L^2(\Om\setminus\omega)}$ converges to $0$ as $k$ goes to infinity, for every $j\in \N$, we have that $\|\sqrt{|W|}\phi_j(\Om,V_k)\|_{L^2(\Om\setminus\omega)}$ converges to $0$ as $k$ goes to infinity and, by equation~\eqref{eq-disug}, we deduce  that there exists $\bar k$ large enough such that $V_{\bar k}\in\hcR_n(\Om,W)$.

By Proposition~\ref{teytel-docet} there exists an analytic
function $\mu\mapsto W_\mu$  from $[0,1]$ into $L^\infty(\Omega)$
such that $W_0=0$, $W
_1=V_{\bar k}-\overline V$ and the spectrum of $-\Delta+\overline V+W_\mu$
is simple for every $\mu\in(0,1)$. Therefore applying Proposition~\ref{sofisticazzi} and since  $V_{\bar k}=\overline V+W_1\in \hcR_n(\Om,W)$
we get that $\overline V+W_\mu\in \hcR_n(\Om,W)$ for almost every $\mu\in (0,1)$, so that  $\hcR_n(\Om,W)$ is dense in $\hat  \calV(\Om,W)$.

The set $\cap_{n\in\N}\hcR_n(\Om,W)$ is then residual in 
$\hat\calV(\Om,W)$
. We claim that if $V\in \cap_{n\in\N}\hcR_n(\Om,W)$ then there exists a reordering $\hat h$ of $\N$ such that 
$B^{\hat h}_{n}(\Om,V,W)$ is connected for every $n\in\N$. 
Indeed, let $\al$ be a map from the power set of $\N$ into itself
defined by
\[
\al(J)=\{m\in\N\setminus J\mid\int_\Om W \phi_{n}(\Om,V)\phi_{m}(\Om,V)\ne0 \mbox{ for some }n\in J\}.
\]
Then $\hat h$ can be defined inductively as follows: set $\hat h(1)=1$ and, for every $n\in\N$,  
let $\hat h(n+1)$ be the smallest element of $\alpha(\{\hat h(1),\dots,\hat h(n)\})$. 
It is straightforward to check that $\hat h$ is a reordering of $\N$.

The triple 
$(\Om,V,W)$
is then \effective\ if $V$ belongs to 
$$\lp \cap_{n\in\N}\hcR_n(\Om,W)\rp\cap \lp \cap_{q\in \cup_{K\in \N}\Q^K\setminus \{0\}}\OO_q(\Om)\rp$$
that is the intersection of countably many open and dense subsets of $\hat\calV(\Om,W)$.
\EOP

\begin{rem}\label{tuttouguale}
The proof shows that if $\sigma(\Om,V)$ is simple and if, 
for  every pair of indices $j,k\in\N$, there exists a finite sequence $r_1,\dots,r_l\in\N$ such that 
$r_1=j$, $r_l=k$,  and 
$$\int_\Om W \phi_{r_i}(\Om,V)\phi_{r_{i+1}}(\Om,V)\ne 0$$
for every $i=1,\dots,l-1$, 
 then there exists a reordering $h$ of $\N$ for which $B^{h}_{n}(\Om,V,W)$ is connected for every $n\in\N$.  
In particular, if $\sigma(\Om,V)$ is simple and for some 
reordering $h$ the matrices  $B^h_{n}(\Om,V,W)$ are frequently connected (as in the definition of \effectiveness),   
then we can assume without loss of generality 
 that  $B^{h}_{n}(\Om,V,W)$ is connected for every $n\in\N$.  
\end{rem}

The next corollary follows immediately from Theorem~\ref{V-s}.
\begin{corol}
Let $\Omega\in\Xi_d$ and $W\in L^\infty(\Om)$ be non-constant and absolutely continuous. 
Then, generically with respect to
$V$ in 
$L^\infty(\Om)$,
the triple 
$(\Om,V,W)$
is \effective.
\end{corol}
In the unbounded case we deduce the following.
\begin{corol}\label{4.5}
Let $\Omega=\R^d$ and $W\in\calW(\R^d)$ be non-constant and absolutely continuous. Assume that $U\subset\R$ has nonempty interior.
Then, generically with respect to
$V$ in 
$\calV(\R^d,W,U)$,
the quadruple 
$(\R^d,V,W,U)$
is \good.
\end{corol}
\proof
Let $u$ belong to the interior of $U$. Assume in particular $[u-\delta,u+\delta]\subset U$. Then, from the definition of  $\calV(\R^d,W,U)$ we have that $V+uW+\delta W$ and $V+uW-\delta W$ are both positive outside a bounded subset $\Om_0$ of $\R^d$. In particular $|W|\leq \frac1\delta |V+uW|$ outside $\Om_0$, while $W$ is bounded on $\Om_0$.
Therefore $V+uW\in\hat\calV(\R^d,W)$ and applying Theorem~\ref{V-s}, we have that the triple $(\R^d,V+uW,W)$ is  \effective, generically with respect to $V\in \mathcal V(\R^d,W,U)$.
\EOP

\subsection{The quadruple $(\Om,V,W,U)$ is generically \good\
with respect to $W$
}

We prove in this section that for a fixed potential $V$,  generically
 with respect to
 $W\in \calW(\Om,V,U)$, 
the quadruple 
$(\Om,V,W,U)$
is \good.
 Notice that $(\Om,V,W)$ cannot be \effective\ if the spectrum of $-\Delta+V$ is resonant, independently of $W$. In this regard the result is necessarily weaker than Theorems~\ref{wrtVW} and \ref{V-s}, where the genericity of the \effectiveness\
 was proved.

\begin{prop}\label{W-s}
Let $\Omega$ belong to $\Xi_d^\infty$ and
$V\in\calV(\Om)$ be an absolutely continuous function on $\Om$.
Assume that $U$ has nonempty interior.
Then, generically with respect to
$W\in\calW(\Om,V,U)$, the quadruple $(\Om,V,W,U)$
is \good.
\end{prop}
\proof
Fix $u\not=0$ in the interior of $U$. 
Notice that $V+u\calW(\Om,V,U)$ is an open subset of $\calV(\Om)$ diffeomorphic to $\calW(\Om,V,U)$.
In particular,
for every $K\in\N$ and $q\in \Q^K\setminus\{0\}$, the set 
$\{W\in \calW(\Om,V,U)\mid V+u W\in\OO_q(\Om)\}$
is open and dense in $\calW(\Om,V,U)$.

For every $W\in \calW(\Om,V,U)$ let $\hcR_n(\Om,W)$ be defined as in the previous section. 
As proved in Corollary~\ref{4.5}, for every $W\in \calW(\Om,V,U)$ 
one has $V+uW\in \hat{\cal V}(\Om,W)$. 
We prove the proposition by showing that
for every  $n\in\N$,
for each $W$  in a open and dense subset of $\calW(\Om,V,U)$ (depending on $n$),
$V+u W$ belongs to $\hcR_n(\Om,W)$.

Define 
$${\cal P}_n=
\{W\in \calW(\Om,V,U)\mid V+u W\in\hcR_n(\Om,W)\}\,.$$
Because of Remark~\ref{tuttouguale} it is enough to prove that each ${\cal P}_n$ is open and dense.
Since 
$$W\mapsto\int_\Om W \phi_j(\Om,V+uW)\phi_k(\Om,V+uW).$$
is continuous on $\{W\in \calW(\Om,V,U)\mid \lb_j(\Om,V+uW),\,\lb_k(\Om,V+uW)\mbox{ are simple}\}$ for every $j,k\in\N$ (Proposition~\ref{the-supercontinuity}),  we deduce that ${\cal P}_n$ is open. 

Fix $\overline{W}\in  \calW(\Om,V,U)$. 
We are left to prove that $\overline{W}$ belongs to the closure of ${\cal P}_{n}$.

Consider first the case in which $V$ is constant. In particular, $\Om\in\Xi_d$, 
$\calW(\Om,V,U)=V+u\calW(\Om,V,U)=L^\infty(\Om)$, and 
\be
\int_\Om W \phi_j(\Om,V+uW)\phi_k(\Om,V+uW)=\int_\Omega W\phi_j(\Om,u W)\phi_k(\Om,u W).
\label{caso-con}
\ee

Fix $\omega\in\Xi_d$ compactly contained in $\Omega$, whose
boundary is Lipschitz continuous and such that the spectrum $\sigma(\omega,0)$ is simple.
For instance, $\om$ can be taken as an orthotope
whose side's lengths are  non-resonant.
(The simplicity of the spectrum of the Laplace-Dirichlet operator on  $\omega$ is
actually generic among sufficiently smooth domains, as
proved in \cite{micheletti-Perturbazione,ule}.)

Let $z\in L^\infty(\om)$ be non-orthogonal in $L^2(\om)$ to
$\phi_j(\om,0)\phi_k(\om,0)$ for every $j,k\in\N$. (Such $z$ exists because
each $\phi_j(\om,0)\phi_k(\om,0)$ is not identically equal to zero and because $L^\infty(\om)$ is a Baire space.)
Then, for every $j,k\in\N$, the derivative of
$$\eps\mapsto \int_\om \eps z\phi_j(\om,\eps z)\phi_k(\om,\eps z)$$
at $\eps=0$ is equal to
$$\int_\om z\phi_j(\om,0)\phi_k(\om,0)\ne 0.$$
By Theorem~\ref{the-analyticity},  there exists $\bar \eps\in\R$ such that the spectrum $\sigma(\om,\bar \eps z)$ is  simple and
$$\int_\om \bar\eps z\phi_j(\om,\bar\eps z)\phi_k(\om,\bar \eps z)\ne0$$
for every $j,k\in\N$.

Let $(W_l)_{l\in\N}$ be a sequence in $L^\infty(\Om)$ such that
$W_l-\overline{W}\in L^\infty(\Om)$, 
$\lim_{l\to\infty}W_l|_\om=(\bar\eps/u) z$ in $L^2(\om)$ and
$\lim_{l\to\infty}\inf_{\Om\setminus\om} W_l=+\infty$.
By Lemma~\ref{buco} we deduce that there exists $\bar l$ large enough such that 
$$\int_\Om W_{\bar l}\phi_j(\Om,uW_{\bar l})\phi_k(\Om,uW_{\bar l})\ne0\quad \mbox{for }j,k=1,\dots,n.$$
By Proposition~\ref{teytel-docet} we can consider an analytic curve $\mu\mapsto \hat W_\mu$  in $L^\infty(\Om)$ for $\mu\in[0,1]$ such that $\hat W_0=\overline{W}$,  $\hat W_1=W_{\bar l}$ and the spectrum of $-\Delta+u\hat W_\mu$ is simple for every $\mu\in (0,1)$, and we have  
$$\int_\Om \hat W_\mu\phi_j(\Om,V+u\hat W_\mu)\phi_k(\Om,V+u\hat W_\mu)=\int_\Om \hat W_\mu\phi_j(\Om,u\hat W_\mu)\phi_k(\Om,u\hat W_\mu)\ne0$$
for almost every $\mu\in(0,1)$ and in particular for some $\mu$ arbitrarily small,  implying that
$\overline{W}$ belongs to the closure of ${\cal P}_{n}$.

Let now $V$ be non-constant.
Let $\om\subset \Om$ and $h$ be as in the statement of Lemma~\ref{AC} with $V$ playing the role of  $Z$.

Take a sequence $(W_k)_{k\in\N}$  in $\calW(\Om,V,U)$ such that $W_k-\overline{W}$ belongs to $L^\infty(\Om)$ for every $k$ and
$$\lim_{k\to+\infty} \|V+uW_k\|_{L^\infty(\om)}=0,\quad \lim_{k\to+\infty}\inf_{\Omega\setminus\omega}(uW_k)=+\infty.$$
According to Lemma~\ref{buco},
\[
\lim_{k\to +\infty}\phi_{m}(\Om,V+uW_k)=\phi_{m}(\om,0)\quad\mbox{and}\quad\lim_{k\to +\infty}\sqrt{V+u W_k}\phi_{m}(\Om,V+uW_k)=0
\]
in $L^2(\Om,\C)$ for every $m\in \N$, where $\phi_{m}(\om,0)$ is identified with its extension by zero on $\Om\setminus\om$. 
In particular, we have that $\sqrt{V}\phi_{m}(\Om,V+uW_k)$ converges in $L^2(\Om,\C)$  as $k$ tends to infinity to the extension by zero of $\sqrt{V}\phi_{m}(\om,0)$ on $\Om\setminus\om$.  
Hence, 
\[
\lim_{k\to +\infty}\int_\Om W_k\phi_{h(l)}(\Om,V+uW_k)\phi_{h(l+1)}(\Om,V+uW_k)=-\frac 1u\int_\om V\phi_{h(l)}(\om,0)\phi_{h(l+1)}(\om,0)\ne0,
\]
for every $l\in\N$. For a fixed $n\in\N$, we can choose $\bar k$ large enough so that
$$\int_\Om W_{\bar k}\phi_{h(l)}(\Om,V+uW_{\bar k})\phi_{h(l+1)}(\Om,V+uW_{\bar k})\neq 0\,,$$ 
for $l$ large enough, in order to guarantee that $W_{\bar k}\in{\cal P}_n$. 
By Proposition~\ref{teytel-docet} there exists an analytic
path $\mu\mapsto \hat{W}_\mu$  from $[0,1]$ into $L^\infty(\Omega)$
such that $\hat{W}_0=0$, $\hat{W}
_1=W_{\bar k}-\overline{W}$ and the spectrum of $-\Delta+V+u \overline W+u \hat{W}_\mu$
is simple for every $\mu\in(0,1)$. Therefore, by analyticity of the eigenfunctions and by applying Proposition~\ref{sofisticazzi}, we get that
$$\int_\Om (\overline W+\hat W_{\mu})\phi_{h(l)}(\Om,V+u \overline W+u\hat W_{\mu})\phi_{h(l+1)}(\Om,V+u \overline W+u\hat W_{\mu})\neq 0$$
for almost every  $\mu\in(0,1)$.

Hence, $\overline{W}$ belongs to the closure of ${\cal P}_n$.
\EOP

\begin{rem} 
It seems possible to adapt 
the arguments presented above and in the previous section 
to the 
conditions ensuring approximate controllability in 
the recent work by Nersesyan 
\cite{nersesyan}: namely, that there exists a reordering $h$ such that
$\lb_{h(1)}(\Om,V)-\lb_{h(j)}(\Om,V)\not=\lb_{h(p)}(\Om,V)-\lb_{h(q)}(\Om,V)$ 
 for all $j,p,q\in\N$ such that $\{1, j \}\not =\{p, q\}$ and $j\not=1$,
 and
that 
the first line of $B^h_n(\Om,V,W)$ is made of non-zero elements for every $n\in \N$ (\cite{nersesyan} also requires that $\Om$ is bounded, with smooth boundary and that $V,W$ are smooth up to the boundary). 
In order to do so, a counterpart of Lemma~\ref{AC} should be proved, replacing
\r{eq_snake} by
$$\int_\om Z\phi_{h(1)}(\om,0)\phi_{h(l)}(\om,0)\ne 0,\quad\mbox{for every $l\in\N$}.$$
This is done in \cite{BCKL} for the case $d=2$ (just replace $\mu$ by $(Z,0)$ in Proposition~2.8). 
 \end{rem}

\section{Conclusion}\label{concl}
In this paper we proved
that once $(\Om,V)$ or $(\Om,W)$ is fixed (with $W$ non-constant), the bilinear Schr\"odinger equation on
$\Om$ having $V$ as uncontrolled and $W$ as controlled potential is generically approximately controllable in the sense of the density matrices
with respect to the other element of the triple $(\Om,V,W)$.

A natural question is whether a similar property holds with respect to the dependence
on $\Om$. It makes sense to conjecture that it does but the proof of this fact seems
hard to obtain through the techniques used here.
Fix $V$ and $W$
absolutely continuous on $\R^d$
with $W$ nowhere locally constant.
Let $m\in\N$ and $\Om$ belong to the space of bounded ${\cal C}^m$ domains endowed with the 
${\cal C}^m$ topology (this space is Baire as proved in \cite{micheletti-baire}).
One important remark is that the dependence of $\lb_k(\Om,V)$ on $\Om$ is not necessarily analytic, as it would be the case if $V$ was analytic. (A genericity non-resonance result for the spectrum in the case
$V=0$, for instance,  has been proved along these lines in \cite{yannick}.)
Similarly, the quantities $\int_\Om W \phi_k(\Om,V)\phi_j(\Om,V)$ do not in general vary analytically with respect to $\Om$.
Hence, the pattern of the proofs seen in the previous sections could not be followed.
A partial result going in the right direction can be found in \cite{BCKL}, where the authors prove that for $V=0$ and $W\in{\cal C}^1(\R^2,\R)$ nowhere-constant,
for a generic ${\cal C}^3$ domain $\Om\subset \R^2$ one has $\int_\Om W \phi_1(\Om,0)\phi_j(\Om,0)\ne 0$ for every $j\in\N$.
The proof of this fact in \cite{BCKL} is very technical and ingenious. Its extension to 
general uncontrolled potentials and to the case $d>2$ looks an extremely hard task.

\section{Acknowledgments}
The authors are grateful to
Ugo Boscain, Thomas Chambrion, Yacine Chitour,  
Antoine Henrot and Karim Ramdani for helpful discussions.

\bibliographystyle{abbrv}
\bibliography{biblio_g}




%

\end{document}